\documentclass{article}

\usepackage[utf8]{inputenc} 
\usepackage[T1]{fontenc}   
\usepackage{amsmath}        
\usepackage{amssymb}      
\usepackage{amsfonts}       
\usepackage{amsthm}         
\usepackage{graphicx}      
\usepackage[margin=1in]{geometry}
\usepackage{hyperref}       
\usepackage{booktabs}      
\usepackage{caption}       
\usepackage{xcolor}        
\usepackage{subcaption}

\newcommand{\ap}{a_p(E)}        

\definecolor{zvoColor}{RGB}{0, 200, 0}

\title{Improving elliptic curve rank classification using multi-value and learned Mestre–Nagao sums}

\author{
  Zvonimir Bujanović\thanks{University of Zagreb, Faculty of Science,
                           Department of Mathematics} \and
  Matija Kazalicki\footnotemark[1] \and
  Domagoj Vlah\thanks{University of Zagreb, Faculty of Electrical Engineering and Computing, Department of Applied Mathematics}
}

\begin{document}
\date{}
\maketitle

\begin{abstract}
Determining the rank of an elliptic curve $E/\mathbb{Q}$ remains a central challenge in number theory. Heuristics such as Mestre--Nagao sums are widely used to estimate ranks, but there is considerable room for improving their predictive power. This paper introduces two novel methods for enhancing rank classification using Mestre--Nagao sums. First, we propose a ``multi-value'' approach that simultaneously uses two distinct sums, $S_0$ and $S_5$, evaluated over multiple ranges. This multi-sum perspective significantly improves classification accuracy over traditional single-sum heuristics. Second, we employ machine learning---specifically deep neural networks---to learn optimal, potentially conductor-dependent weightings for Mestre--Nagao sums directly from data. Our results indicate that adaptively weighted sums offer a slight edge in rank classification over traditional methods.

\end{abstract}

\section{Introduction}
\label{sec:introduction}
The rank \( r \) of the Mordell-Weil group of an elliptic curve \( E \) over the rational numbers \( \mathbb{Q} \) is a fundamental invariant that quantifies the size of the group of rational points \( E(\mathbb{Q}) \). Determining this rank is a central challenge in number theory. While algorithms exist for computing the rank, they often entail significant computational complexity and may not terminate for all curves. The Birch and Swinnerton-Dyer (BSD) conjecture connects the rank to the analytic properties of the curve's $L$-function, inspiring the development of analytic heuristics for rank estimation.

One such heuristic is the Mestre-Nagao sum \( S_0(B) \), defined as a sum over the primes \( p \leq B \) not dividing the conductor \( N \) of the curve:

\[
S_0(B) = \frac{1}{\log B} \sum_{\substack{p \leq B \\ p \nmid N}} \frac{a_p \log p}{p},
\]
where \( a_p = p + 1 - \#E(\mathbb{F}_p) \) is the trace of Frobenius. Under the original BSD conjecture, it is expected that \( \lim_{B \to \infty} S_0(B) = -r_{\text{an}} + \frac{1}{2} \), where \( r_{\text{an}} \) denotes the analytic rank of \( E \), conjecturally equal to the algebraic rank \( r \)(see \cite{KimMurty23}).

In practice, \( S_0(B) \) is computed for a large bound \( B \), and its value is used to predict the rank. However, the effectiveness of this sum can be limited. Its convergence can be slow (see \cite[Appendix]{KimMurty23}), and recent work \cite{BKN24, POZ24} has highlighted an
oscillatory behavior linked to the phenomenon known as \textbf{murmurations}
of elliptic curves \cite{HLOP24, LOP25, Zub23}. This phenomenon surprisingly
reveals that increasing the summation bound $B$ does not always lead to
better rank classification accuracy; in some cases, smaller bounds yield
superior results. This observation underscores the subtleties
in applying these sums naively. Furthermore, the standard Mestre-Nagao sum
treats all primes equally (up to the weighting) and does not adapt to
curve-specific properties like the conductor $N$.

Recent studies have applied machine learning, notably linear regression \cite{HLO23} and Convolutional Neural Networks (CNNs) \cite{KV23}, to the problem of rank classification using \( a_p \) sequences. In a similar spirit, several new works \cite{BCCHLNP25, BCDLLOQV25, BBFHHS24} have explored predicting other deep invariants of elliptic curves—such as Euler factors, orders of vanishing of \( L \)-functions, and properties of the Shafarevich–Tate group—using machine learning approaches. Pursuing related goals, this paper proposes two novel techniques specifically designed to enhance Mestre–Nagao sum heuristics.

First, we propose a \textit{multi-value Mestre--Nagao sum} approach. Instead of relying on a single sum \( S(B) \) computed at a fixed bound \( B \), we consider multiple sums, such as \( S_0(B_i) \) and \( S_5(B_j) \)(for definition see Section \ref{sec:multi_value}), evaluated at various bounds \( B_i \) and \( B_j \).   This collection of sum values serves as input features for a classification model in a neural network, resulting in improved discrimination between different ranks compared to using any single sum value alone (see Table \ref{tab:mcc_top_range}).

Second, we propose a machine learning approach for constructing an \textit{optimal Mestre-Nagao sum}, a weighted combination of Frobenius traces, where the weights may depend on the conductor, depending on whether it is included as an input feature. Specifically, we study expressions of the form
\[
\sum_{p \leq B} w_p \ap/\sqrt{p}
\]
where the weights $w_p$ are learned from data to maximize the accuracy of rank classification. This approach allows the model to adaptively emphasize the contributions of different primes.

We evaluate these approaches using LMFDB \cite{LMFDB} and Balakrishnan et al. \cite{BHKSSW16} database of elliptic curves,
comparing their performance against traditional Mestre-Nagao heuristics implemented
by machine learning models. For implementing our models we use PyTorch \cite{PGMLBC19}. Our results indicate that both the
multi-value approach and the learned sums offer tangible improvements in
rank classification accuracy.

This paper is structured as follows: Section~\ref{sec:multi_value} details
the multi-value Mestre-Nagao sum approach and presents classification
results. Section~\ref{sec:learning_optimal} describes the methodology for
learning optimal sum coefficients using neural networks and analyzes the
learned weights and resulting performance. Section~\ref{sec:conclusion}
concludes with a summary of findings and potential directions for future
research.

\section{Rank Classification via Multi-Value Mestre-Nagao Sums}
\label{sec:multi_value}

In this section, we focus on the Mestre-Nagao sums $S_0(B)$
and $S_5(B)$, defined as:
\begin{align}
    S_0(B) &=\frac{1}{\log{B}} \sum_{\substack{p \le B \\ p \nmid N}} \frac{\ap \log p}{p}
             \label{eq:unnormalized_mestre_nagao_0}, \\
    S_5(B) &= \sum_{\substack{p \le B \\ p \nmid N}} \log\left(\frac{p+1-\ap}{p}\right)+\sum_{\substack{p \le B \\ p | N}}\log\left(1.5\,\frac{p - 1}{p}\right)
        \label{eq:unnormalized_mestre_nagao_5}.
\end{align}
A variant of $S_5$ was employed in \cite{EK20} to find elliptic curves with record ranks.

We investigate the use of multiple Mestre--Nagao sums, specifically \( S_0(B) \) and \( S_5(B) \), computed at various bounds \( B \in \{1000, 5000, 10000, 20000, 30000, 40000, 50000, 100000\} \), as input features for rank classification models (together with conductor $N$). These sums are combined in different configurations to capture a richer set of information about the curve's behavior. Fully connected neural networks are employed as the classification models.

More precisely, consider $\mathcal{B} \subseteq \{1000, 5000, 10000, 20000, 30000, 40000, 50000, 100000\}$, and a neural network that takes as input the logarithm of the conductor $N$ and the Mestre-Nagao sums $S_0(B)$ and $S_5(B)$ for all $B \in \mathcal{B}$. Several fully connected layers with ReLU activations are used as the hidden layers of the network, where the optimal number of layers and the number of neurons in each hidden layer were determined by hyperparameter search. We used ranges of $3$ to $6$ hidden layers and $8$ to $512$ neurons in each hidden layer (using only powers of $2$). Best results are achieved by using $4$ hidden layers and $64$ or $128$ hidden neurons per layer for experiments in subsection \ref{subsection:TopRangeTest} or $256$ hidden neurons per layer for experiments in subsection \ref{subsection:UniformTestRange}. The output layer computes the classification probabilities for each of the possible ranks. A weighted cross-entropy loss function is used in the optimization of this neural network, where weights are computed to be
proportional to the inverse of the relative frequency of each rank in the used elliptic curve dataset. The AdamW optimizer \cite{LoshchilovHutter19} was used to train the network, with a hyperparameter search to find the best learning rate and weight decay. We also consider the case where only the $S_0$ sums are used as input, and the case where only the $S_5$ sums are used as input (together with the logarithm of the conductor in all cases).

Two primary experimental setups were used, as described in the following subsections.

\subsection{Top Range Test}\label{subsection:TopRangeTest}
In the first setup, models were trained on elliptic curves from the 
Balakrishnan et al. \cite{BHKSSW16} database
with conductors in the range $N \in [1, 10^8]$
(4,875,676 curves). These trained models were then tested on curves with
conductors in the higher range $N \in (10^8, 10^9]$ (12,512,753 curves).

The performance, measured by the Matthews Correlation Coefficient (MCC), for
different combinations of input sums and bounds $B$ is presented in
Table~\ref{tab:mcc_top_range}. The results suggest that using combinations
of sums, particularly $S_0$ and $S_5$ together, generally yields better
performance than using either sum alone. Notably, using all available bounds
for both $S_0$ and $S_5$ achieved the highest MCC (0.795), although using
just $B \in \{10^3, 10^5\}$ still provided strong results (MCC=0.784).
Using only $S_0(10^5)$ gave an MCC of 0.686.

\begin{table}[htbp]
    \centering
    \caption{MCC for neural networks trained with different combinations of sums
             (Top Range Test). Tested with curves in the conductor range
             $N \in (10^8, 10^9]$.}
    \label{tab:mcc_top_range}
    \begin{tabular}{c c c c}
        \toprule
        $ B $                 & $ S_0 $ and $ S_5 $ & $ S_0 $ & $ S_5 $ \\
        \midrule
        All                   & 0.795               & 0.737   & 0.608   \\
        $ 10^3, 10^4, 10^5 $  & 0.786               & 0.726   & 0.609   \\
        $ 10^3, 10^5 $        & \textbf{0.784}      & 0.709   & \textbf{0.595}   \\
        $ 10^4, 10^5 $        & 0.740               & 0.705   & 0.562   \\
        $ 5 \cdot 10^4, 10^5 $& 0.710               & 0.685   & 0.537   \\
        $ 10^5 $              & 0.699               & \textbf{0.686}   & \textbf{0.527}   \\
        $ 10^4 $              & 0.553               & 0.537   & 0.384   \\
        $ 10^3 $              & 0.330               & 0.309   & 0.277   \\
        \bottomrule
    \end{tabular}
\end{table}

Confusion matrices provide further insight into the classification performance
for specific rank predictions. Table~\ref{tab:cm_top_range_s0s5allB} shows the
confusion matrix (in percentages) for the best performing model using $S_0$,
$S_5$, and all values of $B$. Table~\ref{tab:cm_top_range_s0_1k} shows the
matrix for a baseline model using only $S_0(1000)$.

\begin{table}[htbp]
    \centering
    \caption{Confusion matrix (\%) for neural network predictions with $ S_0 $,
             $ S_5 $, and all values of $ B $ (Top Range Test).
          }
    \label{tab:cm_top_range_s0s5allB}
    \begin{tabular}{c c c c c c}
        \toprule
        True Rank & Pred 0 & Pred 1 & Pred 2 & Pred 3 & Pred 4 \\
        \midrule
        0 & 23.674 &  6.497 &  0.336 & 0.000 & 0.000 \\
        1 &  3.854 & 41.410 &  1.578 & 0.003 & 0.003 \\
        2 &  0.474 &  0.521 & 18.352 & 0.003 & 0.000 \\
        3 &  0.000 &  0.001 &  0.005 & 3.134 & 0.000 \\
        4 &  0.000 &  0.000 &  0.000 & 0.000 & 0.153 \\
        \bottomrule
    \end{tabular}
\end{table}

\begin{table}[htbp]
    \centering
    \caption{Confusion matrix (\%) for neural network predictions with $ S_0 $
             only, $ B = 1000 $ (Top Range Test).}
           
    \label{tab:cm_top_range_s0_1k}
    \begin{tabular}{c c c c c c}
        \toprule
        True Rank & Pred 0 & Pred 1 & Pred 2 & Pred 3 & Pred 4 \\
        \midrule
        0 & 6.607 & 17.907 &  5.977 & 0.016 & 0.000 \\
        1 & 4.403 & 27.587 & 14.301 & 0.558 & 0.000 \\
        2 & 0.000 &  2.602 & 15.764 & 0.984 & 0.000 \\
        3 & 0.000 &  0.000 &  0.462 & 2.677 & 0.000 \\
        4 & 0.000 &  0.000 &  0.000 & 0.043 & 0.110 \\    
        \bottomrule
    \end{tabular}
\end{table}

\subsection{Uniform Test Range}\label{subsection:UniformTestRange}
In the second setup, models were trained, validated, and tested on elliptic
curves with conductors across the entire range $N \in [1, 10^9]$. The full
dataset of 17,388,429 curves was split into 60\% for training, 20\% for
validation, and 20\% for testing.

The MCC results for this uniform range setup are shown in
Table~\ref{tab:mcc_uniform_range}. Similar trends are observed, with the
combination of $S_0$ and $S_5$ using all bounds achieving the highest MCC
(0.856). Using only $S_0(10^5)$ yielded an MCC of 0.712.

\begin{table}[htbp]
    \centering
    \caption{MCC for neural networks trained with different combinations of sums
             (Uniform Test Range).}
    \label{tab:mcc_uniform_range}
    \begin{tabular}{c c c c}
        \toprule
        $ B $                 & $ S_0 $ and $ S_5 $ & $ S_0 $ & $ S_5 $ \\
        \midrule
        All                   & 0.856               & 0.787   & 0.681   \\
        $ 10^3, 10^4, 10^5 $  & 0.855               & 0.769   & 0.661   \\
        $ 10^3, 10^5 $        & \textbf{0.835}      & 0.747   & \textbf{0.637}   \\
        $ 10^4, 10^5 $        & 0.787               & 0.739   & 0.601   \\
        $ 5 \cdot 10^4, 10^5 $& 0.752               & 0.729   & 0.578   \\
        $ 10^5 $              & 0.733               & \textbf{0.712}   & \textbf{0.573}   \\
        $ 10^4 $              & 0.596               & 0.573   & 0.422   \\
        $ 10^3 $              & 0.371               & 0.359   & 0.295   \\
        \bottomrule
    \end{tabular}
\end{table}

Corresponding confusion matrices for the uniform range test are shown in
Table~\ref{tab:cm_uniform_range_s0s5allB} (for $S_0$, $S_5$, all $B$) and
Table~\ref{tab:cm_uniform_range_s0_1k} (for $S_0$ only, $B=1000$).

\begin{table}[htbp]
    \centering
    \caption{Confusion matrix (\%) for neural network predictions with $ S_0 $,
             $ S_5 $, and all values of $ B $ (Uniform Test Range).
             }
    \label{tab:cm_uniform_range_s0s5allB}
    \begin{tabular}{c c c c c c}
        \toprule
        True Rank & Pred 0 & Pred 1 & Pred 2 & Pred 3 & Pred 4 \\
        \midrule
        0 & 26.035 &  4.104 &  0.345 & 0.000 & 0.000 \\
        1 &  3.443 & 42.346 &  1.237 & 0.002 & 0.000 \\
        2 &  0.061 &  0.147 & 19.148 & 0.002 & 0.000 \\
        3 &  0.000 &  0.000 &  0.000 & 2.998 & 0.000 \\
        4 &  0.000 &  0.000 &  0.000 & 0.000 & 0.131 \\
        \bottomrule
    \end{tabular}
\end{table}

\begin{table}[htbp]
    \centering
    \caption{Confusion matrix (\%) for neural network predictions with $ S_0 $
             only and $ B = 1000 $ (Uniform Test Range).
            }
    \label{tab:cm_uniform_range_s0_1k}
    \begin{tabular}{c c c c c c}
        \toprule
        True Rank & Pred 0 & Pred 1 & Pred 2 & Pred 3 & Pred 4 \\
        \midrule
        0 & 7.329 & 17.556 &  5.559 & 0.040 & 0.000 \\
        1 & 3.463 & 30.043 & 12.592 & 0.931 & 0.000 \\ 
        2 & 0.000 &  2.054 & 15.894 & 1.404 & 0.007 \\
        3 & 0.000 &  0.000 &  0.126 & 2.846 & 0.027 \\
        4 & 0.000 &  0.000 &  0.000 & 0.001 & 0.130 \\
        \bottomrule
    \end{tabular}
\end{table}

Overall, our findings demonstrate that combining \( S_0(B) \) and \( S_5(B) \) across various bounds \( B \) enhances rank classification performance compared to using single sum values.

Tables~\ref{tab:mcc_top_range} and~\ref{tab:mcc_uniform_range} show a notable improvement in classification accuracy when transitioning from a model that uses only \( S_5(10^5) \) as input (MCC scores of 0.527 and 0.573) to one that uses both \( S_5(10^3) \) and \( S_5(10^5) \) (MCC scores of 0.595 and 0.637). This can be heuristically explained by recalling the original Birch and Swinnerton-Dyer conjecture, which states that
\[
\prod_{\substack{p \le B \\ p \nmid N}} \frac{p+1 - a_p}{p} \sim C_E (\log B)^r \quad \text{as } B \to \infty,
\]
for some constant \( C_E \) depending on the curve. This suggests the approximation
\[
S_5(B) \approx \log C_E + r \log \log B.
\]
Thus, providing the network with two values of \( S_5 \) at different bounds may allow it to implicitly estimate the constant \( C_E \), thereby reducing misclassification errors.

\subsection{Rank Prediction via Rectangular Regions in Sum Space}

The benefit of using multiple bounds for the same type of sum, such as the
pair $(S_0(10^3), S_0(10^5))$, can be further analyzed. One can attempt to define classification rules based directly
on these two values. Simple rules, like defining rectangular regions in the
$(S_0(10^3), S_0(10^5))$ plane for each rank, can already outperform using only $S_0(10^5)$
in certain conductor ranges.

\begin{figure}
    \centering
    \begin{subfigure}{0.49\textwidth}
        \centering
        \includegraphics[width=\linewidth]{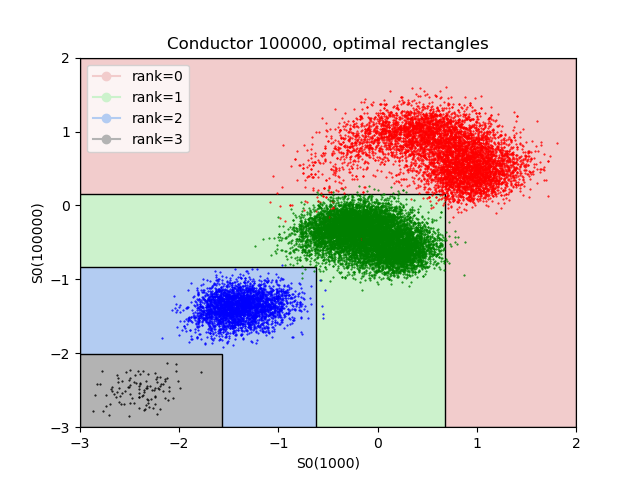}
    \end{subfigure}
    \begin{subfigure}{0.49\textwidth}
        \centering
        \includegraphics[width=\linewidth]{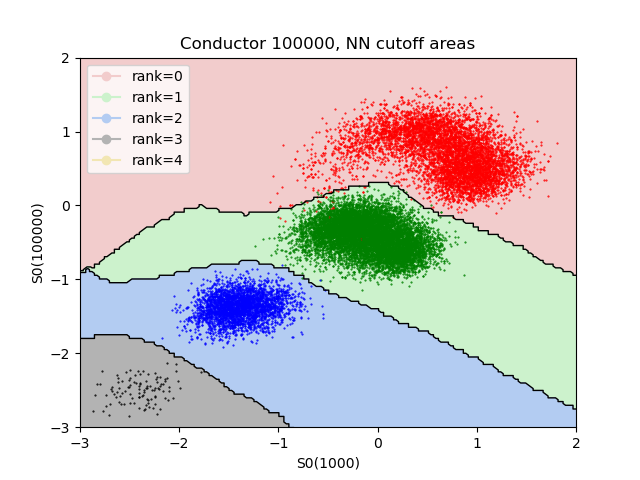}
    \end{subfigure}
    \caption{Rank cutoff areas for elliptic curves of conductor $100\,000$ using the values of $S_0(1\,000)$ and $S_0(100\,000)$. 
        Left: Classification based on optimal rectangular regions. Right: Classification by the neural network model. 
        Each point in the point cloud represents a single elliptic curve whose actual rank is color-coded.}
    \label{fig:twosums-s0-1}
\end{figure}

\begin{figure}
    \centering
    \begin{subfigure}{0.49\textwidth}
        \includegraphics[width=\linewidth]{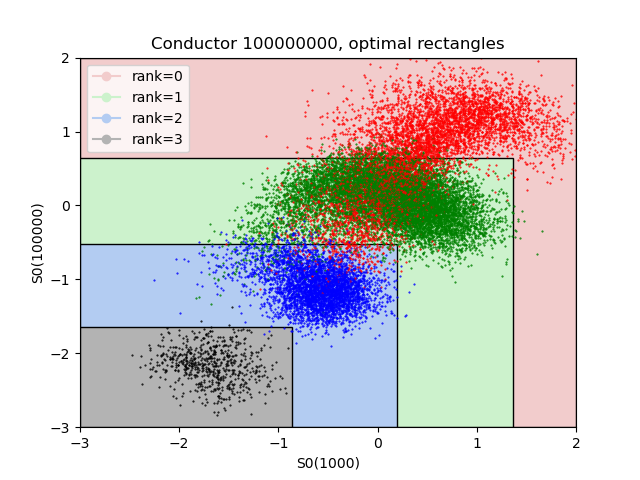}
    \end{subfigure} \hfill
    \begin{subfigure}{0.49\textwidth}
        \includegraphics[width=\linewidth]{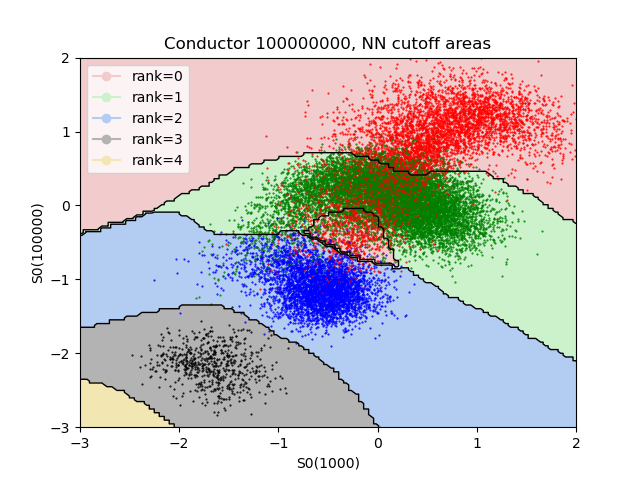}
    \end{subfigure}
    \caption{Rank cutoff areas for elliptic curves of conductor $100\,000\,000$ using the values of $S_0(1\,000)$ and $S_0(100\,000)$. 
        Left: Classification based on optimal rectangular regions. Right: Classification by the neural network model.
        Each point in the point cloud represents a single elliptic curve whose actual rank is color-coded.}
    \label{fig:twosums-s0-2}
\end{figure}

This is demonstrated in Figures~\ref{fig:twosums-s0-1} (left),
and \ref{fig:twosums-s0-2} (left). In these figures, we collected all elliptic curves with conductor $N \in [100\,000, 150\,000]$ (Figure~\ref{fig:twosums-s0-1}) and $N \in [100\,000\,000, 100\,750\,000]$ (Figure~\ref{fig:twosums-s0-2}), and plotted the values of $S_0(10^3)$ and $S_0(10^5)$ for each curve. The color of each point indicates the actual rank of the curve. We then defined rectangular regions in the $(S_0(10^3), S_0(10^5))$ plane for each rank such that the MCC of the classification based on these regions is maximized. The resulting classification regions are shown in the left panels of Figures~\ref{fig:twosums-s0-1} and \ref{fig:twosums-s0-2}. Compared to using only $S_0(10^5)$, this approach improves the MCC from 0.9826 to 0.9921 for the first example (Figure~\ref{fig:twosums-s0-1}) and from 0.7530 to 0.7649 for the second example (Figure~\ref{fig:twosums-s0-2}). Note that the rectangles have to be recomputed for each conductor range.

On the other hand, the neural network trained on the logarithm of the conductor together with $S_0(10^3)$ and $S_0(10^5)$ can learn more complex decision boundaries than simple rectangles, as shown in the right panels of Figures~\ref{fig:twosums-s0-1} and \ref{fig:twosums-s0-2}. The model reported in Table~\ref{tab:mcc_uniform_range} is used for both figures, and the MCC for the elliptic curves from the stated conductor ranges has now improved to 0.9962 for the first example (Figure~\ref{fig:twosums-s0-1}) and to 0.7745 for the second example (Figure~\ref{fig:twosums-s0-2}). Note the intricate zone of rank 0 curves appearing in between the rank 1 and rank 2 zones in the second example which is captured by the neural network model. However, there are many rank 0 curves inside the rank 1 zone, which the neural network model does not detect. This is a known issue with the Mestre-Nagao sums, as they do not always distinguish well between rank 0 and rank 1 curves.

\section{Learning Optimal Mestre-Nagao Sums}
\label{sec:learning_optimal}

While the multi-value approach exploits existing sum definitions more
effectively, it does not alter the underlying structure of the Mestre-Nagao
sums themselves (e.g., the $\log p / p$ weighting). An alternative approach
is to learn optimal weights directly from data.

Consider a generalized weighted sum related to the Mestre-Nagao sum:
\begin{equation} \label{eq:optimal_sum}
S_{opt}(E) = \sum_{p \le B} w_p \frac{\ap}{\sqrt{p}}
\end{equation}
where the weights $w_p$ are parameters learned by the network.
The goal is to find weights $w_p$ such that $S_{opt}(E)$ is maximally
correlated with, or predictive of, the rank of $E$.

We propose using neural networks to learn these optimal weights implicitly.
The network takes as input the sequence of normalized Frobenius traces
$(\ap/\sqrt{p})_{p < 10^5}$ for an elliptic curve $E$, potentially along with
the logarithm of the conductor $\log_{10}(N)$, and is trained to predict the
rank $r \in \{0, 1, 2, 3, 4, 5\}$.

We consider two main architectures:
\begin{enumerate}
    \item \textbf{Conductor-Independent Network:} Takes only the sequence
          $(\ap/\sqrt{p})_{p < 10^5}$ as input.
    \item \textbf{Conductor-Dependent Network:} Takes both
          $(\ap/\sqrt{p})_{p < 10^5}$ and $\log_{10}(N)$ as input. This allows
          the network to potentially learn weights $w_p$ that implicitly
          depend on $N$.
\end{enumerate}
The core of the architecture involves a module that computes the weights $w_p$.
This module uses a fixed positional encoding representing the primes $p < 10^5$. Each positional encoding is a real number from $[-1,1]$ computed as $-1+2(\pi(p)/\pi(B))$, where $\pi(x):=|\{p\leq x: p \ \mathrm{is\ prime}\}|$ is the prime‑counting function.
In the conductor-dependent case, $\log_{10}(N)$ (repeated for each prime) is
concatenated with the positional encoding. This combined input (or just the
positional encoding in the conductor-independent case) is processed through
a series of five 1D convolutional neural network layers (with kernel size 1, acting as per-position
linear transformations) and ReLU activation functions in between layers to produce the weights $w_p$. The first convolutional layer increases the number of channels from $1$ (or $2$ in the conductor-dependent case) to $128$, subsequent convolutional layers keep $128$ channels, and the last convolutional layer decreases the number of channels to $1$, which is the value $w_p$. The generalized sum $S_{opt}(E) = \sum_{p < 10^5} w_p (\ap/\sqrt{p})$ is then
computed. Finally, the pair $(\log_{10}(N), S_{opt}(E))$ is passed through the second module that is composed of four fully connected (Dense) layers with ReLU activations in between layers and $128$ neurons in each of the hidden layers to produce the
final classification probabilities for each of the possible ranks. A weighted cross-entropy loss function is used in the optimization of this neural network, where weights are computed to be proportional to the inverse of the relative frequency of each rank in the used elliptic curve dataset. Again, the AdamW optimizer \cite{LoshchilovHutter19} was used to train the network, with a learning rate equal to $10^{-4}$ and training for $5$ epochs using One Cycle Policy \cite{SmithTopin18} for learning rate scheduling.

\begin{figure}[htbp]
    \centering
    \includegraphics[width=\textwidth]{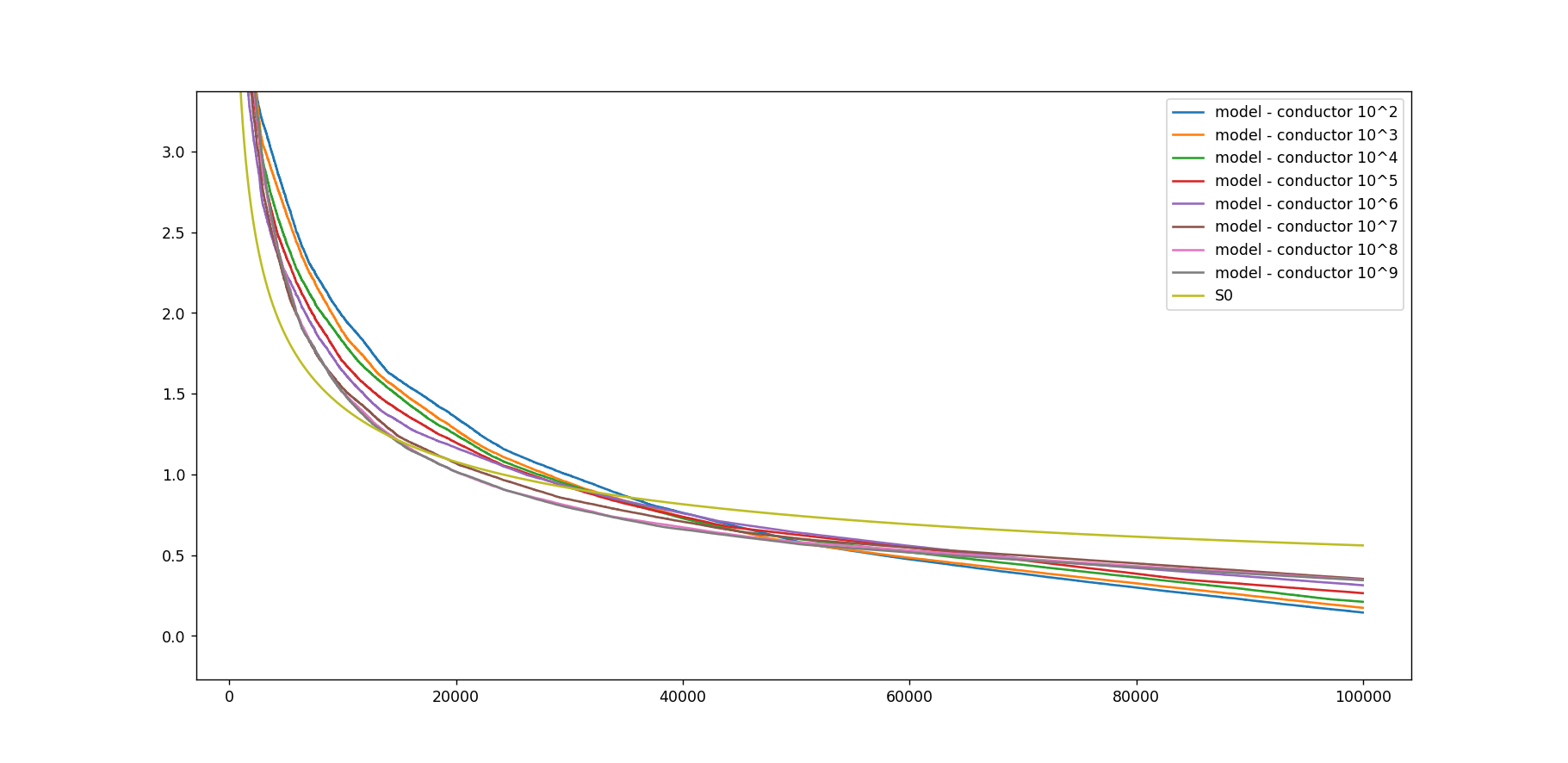}
    
    \caption{Conceptual visualization of learned optimal coefficients $w_p$
             showing dependence on conductor range $N$.}
    \label{fig:learned_coeffs_dependent}
\end{figure}

A key advantage of the conductor-dependent approach is that it allows us to
observe how the implicitly learned optimal weights $w_p$ vary depending on the
conductor range (as illustrated in
Figure~\ref{fig:learned_coeffs_dependent}), providing insight into the
model's learned strategy.

\begin{figure}[htbp]
    \centering
    \includegraphics[width=\textwidth]{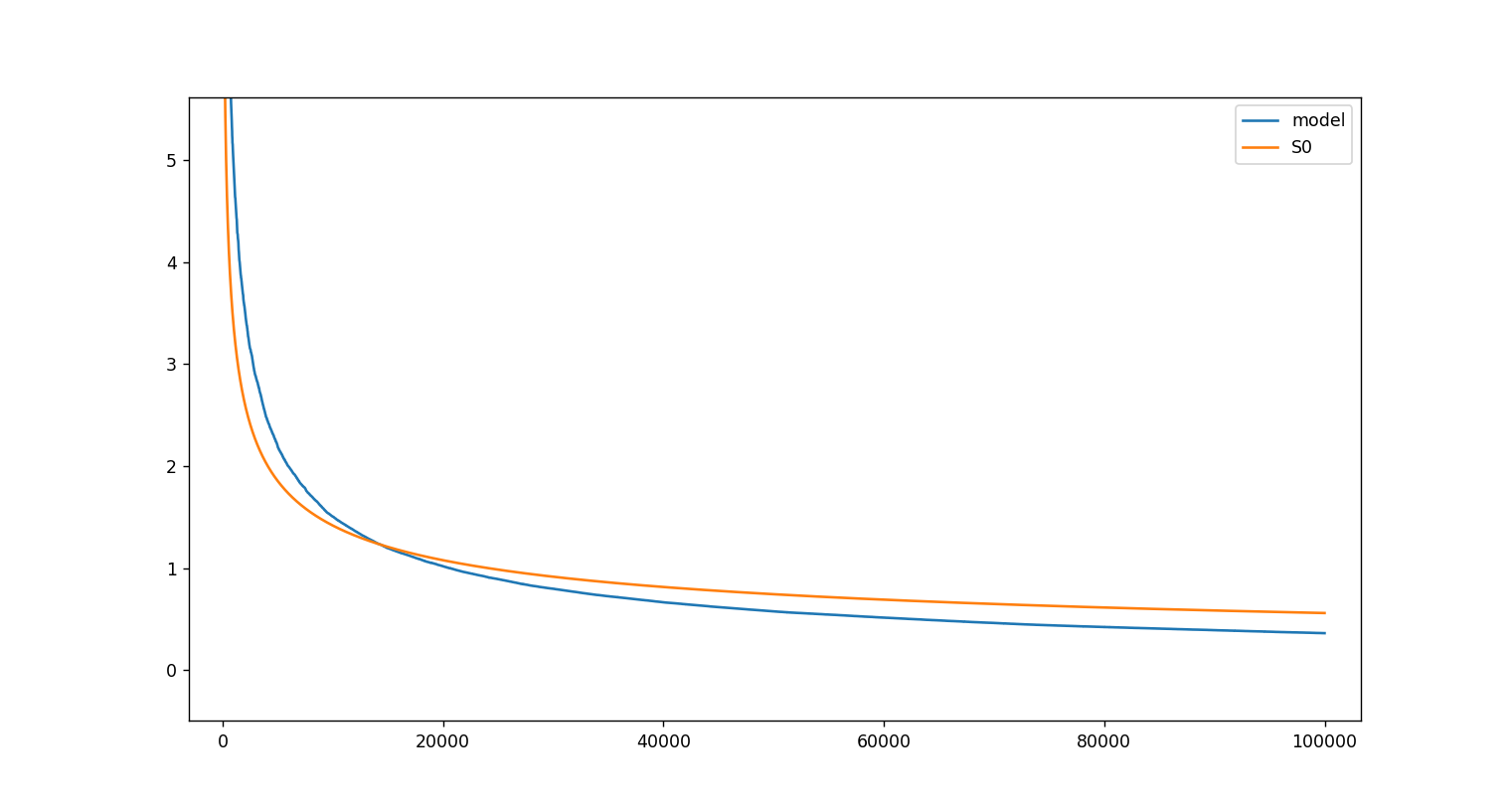}
    
    \caption{Conceptual visualization of learned optimal coefficients $w_p$
             showing dependence on conductor range $N$.}
    \label{fig:learned_coeffs_independent}
\end{figure}

The performance of these learned sum networks was evaluated on the full dataset
($N \in [1, 10^9]$ split 60/20/20), using a uniformly chosen test set.
The conductor-independent network achieved an MCC of $0.7283$, which is the small improvement over the MCC of $0.712$ obtained from classification based on the $S_0$ sum. For the learned coefficients $w_p$ see Figure~\ref{fig:learned_coeffs_independent} The
conductor-dependent network achieved a slightly higher MCC of $0.7322$.

The detailed confusion matrix for the conductor-dependent network
is shown in Table~\ref{tab:confusion_matrix_learned}.

\begin{table}[htbp]
    \centering
    \caption{Confusion matrix (\%) for the rank classification task using the
             learned optimal sums (conductor-dependent network). Ranks 0-5.}
    \label{tab:confusion_matrix_learned}
    \begin{tabular}{c c c c c c c}
        \toprule
        True Rank & Pred 0 & Pred 1 & Pred 2 & Pred 3 & Pred 4 & Pred 5 \\
        \midrule
       
        0 & 29.2371 & 0.9427 & 0.0000 & 0.0000 & 0.0000 & 0.0000 \\
        1 & 4.9485  & 39.8230& 2.1971 & 0.0029 & 0.0000 & 0.0000 \\
        2 & 0.1512  & 0.0007 & 19.1919& 0.0089 & 0.0000 & 0.0000 \\
        3 & 0.0000  & 0.0000 & 0.0001 & 3.0194 & 0.0000 & 0.0000 \\
        4 & 0.0000  & 0.0000 & 0.0000 & 0.0000 & 0.1299 & 0.0000 \\
        5 & 0.0000  & 0.0000 & 0.0000 & 0.0000 & 0.0000 & 0.0009 \\
        \bottomrule
    \end{tabular}
\end{table}

While identifying explicit analytical formulas for these optimal,
conductor-dependent coefficients remains challenging, the neural network
approach effectively discovers and utilizes them. The model's performance on
higher ranks (3, 4, 5), as seen in the confusion matrix,
is limited by the scarcity of such curves in the training data but
demonstrates the potential applicability.

The experiments indicate that neural networks can successfully learn
effective, adaptive weighting schemes for $\ap$ traces, leading to improved
rank classification compared to both traditional Mestre-Nagao sums and the
multi-value approach.

\section{Conclusion and Future Work}
\label{sec:conclusion}
This paper presents two data-driven approaches for improving the rank classification of elliptic curves \( E/\mathbb{Q} \) based on information derived from Frobenius traces \( a_p \). Both methods extend the classical Mestre–Nagao sum heuristic by integrating modern machine learning techniques to enhance predictive accuracy.

The first approach, based on combining multiple Mestre–Nagao sums computed at different bounds, demonstrates that using multiple values—especially from both \( S_0(B) \) and \( S_5(B) \)—significantly improves classification performance compared to using a single sum.

The second approach uses deep neural networks to learn optimal linear combinations of Frobenius traces of the form \( \sum w_p a_p/\sqrt{p} \), resulting in a data-driven refinement of the classical Mestre–Nagao sum. We trained two versions of this model: one that includes the conductor \( N \) as an input feature, and one that does not. Both variants outperformed the baseline model based solely on the \( S_0 \) sum, though the improvement in accuracy was modest.

In future work, we plan to investigate the efficacy of these models on elliptic curves of higher rank, as the main challenge in the present study lies in distinguishing curves of rank zero from those of rank one.


\end{document}